\documentstyle[12pt]{article}

\def\beq{\begin{equation}}
\def\ee{\end{equation}}
\def\bib#1{[{\ref{#1}}]}
\def\eps{\varepsilon}
\def\zba{\overline{z}}

\def\bff{\mbox{\boldmath $f$}}
\def\bfnu{\mbox{\boldmath $\nu$}}
\def\bfZ{\mbox{\boldmath $Z$}}
\newtheorem{theorem}{Theorem}[section]
\newtheorem{proposition}{Proposition}[section]
\newtheorem{corollary}{Corollary}[section]
\newtheorem{definition}{Definition}[section]
\newtheorem{lemma}{Lemma}[section]
\begin{document}           


 \begin{titlepage}

	    \title{
Projective generalizations of Lelieuvre's formula}
 \author{B.G. Konopelchenko$^*$ and U. Pinkall$^{**}$ \\
{\em $^*$ Dipartimento di Fisica, Universit\'{a} di Lecce, 
73100 Lecce, Italy}\\ and \\
{\em $^{**}$ Fachbereich Mathematik, TU Berlin} \\
{\em Str. des 17 Juni, 136, D-10623, Berlin, Germany} }

 \maketitle

      \begin{abstract}
Generalizations of the classical affine Lelieuvre formula to surfaces in 
projective three-dimensional space and to hypersurfaces in multidimensional 
projective space are given. A discrete version of the projective Lelieuvre 
formula is presented too.
      \end{abstract}



 \end{titlepage}

\section{Introduction}
\setcounter{equation}{0}
The classical Lelieuvre formula \bib{r1} of affine geometry provides us
a way to construct a surface via the affine conormal vector (see {\em e.g.}
\bib{r2}-\bib{r4}). Namely, it is the relation 
\beq
{\bff}_{\xi}=\sigma {\bfnu} \wedge 
{\bfnu}_{\xi} \;\; , \;\; 
{\bff}_{\eta}=-\sigma {\bfnu} \wedge 
{\bfnu}_{\eta}
\label{eq:1.1}
\ee
between the coordinates ${\bff}$ of a surface in $R^3$ and its 
conormal ${\bfnu}$. The conormal ${\bfnu}$ 
obeys the equation
${\bfnu}_{\xi \eta} \parallel {\bfnu}$ 
and the corresponding Blaschke metric is
$\Omega=2 det |{\bfnu}, {\bfnu}_{\xi},
{\bfnu}_{\eta} | d\xi \, d\eta$.
For an indefinite metric $\sigma =1$ and $\xi$, $\eta$ are 
real-valued asymptotic coordinates while for a  positive-definite metric 
$\sigma=\sqrt{-1}$ and $\xi$ and $\eta$ are complex conjugate to each other: 
$\eta=\overline{\xi}$. The Lelieuvre formula (\ref{eq:1.1}) is an 
effective tool to study surfaces in affine geometry \bib{r2}-\bib{r4}. It's 
generalization to hypersurfaces in $R^{n+1}$ has been given in \bib{r5}. 
The Lelieuvre formula (\ref{eq:1.1}) provides us also a way to define 
integrable deformations of affine surface via the Nizhik-Veselov-Novikov (
NVN) equation \bib{r6}. 

In this paper we present a projective analog of the Lelieuvre's formula.
It is given by
\beq
f \wedge f_{\xi} = \sigma \star \left( \nu \wedge \nu_\xi \right)
\;\; , \;\;
f \wedge f_{\eta} = -\sigma \star \left( \nu \wedge \nu_\eta \right)
\label{eq:1.2}
\ee
where $f \subset P^3$, $\nu \subset P_3$ and $P_3$ is a projective 
space dual to $P^3$, $\star$ denotes 
the Hodge star operation and $\sigma^2=\pm1$. The relation (\ref{eq:1.2}) 
provides an explicit formula (2.29) for $f$ via $\nu$. The projective 
Lelieuvre map (PLM) (\ref{eq:1.2}) manifest also a symmetry between 
$f$ and $\nu$ (projective duality). We derive the compatibility condition 
for (\ref{eq:1.2}) and prove an invariance of full determinants under the 
correspondence (\ref{eq:1.2}). It is shown that (\ref{eq:1.2}) sets up
correspondence between the normalizations of homogeneous coordinates for a
surface  in $P^3$ and its dual surface in $P_3$.  

A PLM for hypersurfaces in  $(n+1)$-dimensional projective space is 
also presented. A discrete version of the PLM for discrete surfaces in 
$P^3$ is given. We present discrete analogs of the projective Fubini forms. 
An affine reduction ({\em i.e.} the corresponding formulae in the gauge 
$f = \left( {\bff}, -1 \right)$) is considered.
For the discrete case we obtain analogs of the Blaschke and 
affine Fubini cubic forms. 

The paper is organized as follows. In section $2$ and $3$ the PLM
(\ref{eq:1.2}) in asymptotic and conjugate line coordinates  is presented 
and studied. The PLM type formulae for hypersurfaces in $P^{n+1}$ are given 
in section $4$. The discrete analog of the PLM (\ref{eq:1.2}) is presented 
in section $5$. In section $6$ we consider an affine "reduction" of the formulae 
derived.

\section{ The projective Lelieuvre map in asymptotic coordinates.}
\setcounter{equation}{0}

Let $P^n$ and $P_n$ be 
projective spaces dual to each other with homogeneous coordinates 
 $f=(f_1, f_2, f_3, \ldots, f_{n+1})$
and  $\nu=(\nu_1, \nu_2, \nu_3, \ldots, \nu_{n+1})$ respectively.
The pairing of $P_n$ and $P^n$ is defined in a standard way:
\beq
\langle f , \nu \rangle =0 \;\; .
\label{eq:2.1}
\ee
We denote $\eps_{i_1 \ldots i_{n+1}}$ an alternating tensor in $P^n$:
$\eps_{12\ldots n+1}=1$. We will denote the wedge product of $m$
vectors as $a_1 \wedge a_2 \wedge \ldots 
\wedge a_m$.
In a fixed basis one has:
\beq
\begin{array}{l}
\left( a_1 \wedge a_2 \wedge \ldots \wedge a_n \right)_i 
\stackrel{def}{=} \left[ a_1, \dots, a_n \right]_i =
\eps_{i i_2 \ldots i_{n+1}} a_{1_{i_2}} \ldots a_{n_{i_{n+1}}} \; ,\\
a_1 \wedge a_2 \wedge \ldots \wedge a_{n+1} 
=det|a_1, a_2, \dots, a_{n+1}| =
\eps_{i_1 i_2 \dots i_{n+1}} a_{1_{i_1}} a_{i_{i_2}} 
\ldots a_{i_{i_{n+1}}} \;\;.
\end{array}
\label{eq:2.2}
\ee
Note that 
\beq
\langle b ,
[a_1, \dots, a_n ] \rangle = det| b, a_1, a_2, \ldots, a_n |
\;\; .
\label{eq:2.3}
\ee
The Hodge star operation $\star$ on skewsymmetric tensor 
$T_{i_1 \ldots i_{\alpha}}$ is defined as usual
\beq
\left( \star T \right) _{i_{k+1} \ldots i_{n+1}} =\frac{1}{k!}
\eps_{i_1 i_2 \ldots i_{n+1}} T_{i_1 i_2 \ldots i_{k}}
\label{eq:2.4}
\ee
where summation over repeated indices is assumed (here and below). 
One has $\left( \star T \right)=(-1)^{k(n-k)} T$.

In this section we consider the three-dimensional case ($n=3$).
Let $f=f(x,y)$, $\nu=\nu(x,y)$
where $x$ and $y$ are real-valued independent variables. 

\begin{definition}
The projective Lelieuvre map PLM $L: \nu \rightarrow f$
is defined by the equations
\beq
f \wedge f_x = \star \left( \nu \wedge \nu_x \right) \;\; , 
\;\;
f \wedge f_y = -\star \left( \nu \wedge \nu_y \right) \;\; .
\label{eq:2.5}
\ee
\end{definition}
The relations (\ref{eq:2.5}) are manifestly invariant under projective 
transformations in $P^3$ and $P_3$.
\begin{proposition}
The inverse PLM $L^{-1}: f \rightarrow \nu$ is given by equation
\beq
\nu \wedge \nu_x = \star \left( f \wedge f_x \right) \;\; , \;\;
\nu \wedge \nu_x = -\star \left( f \wedge f_y \right) \;\; .
\label{eq:2.6}
\ee
\end{proposition}
The formula (\ref{eq:2.6}) is an obvious consequence of (\ref{eq:2.4})
with $n=4$, $k=2$. The formulae (\ref{eq:2.5}), (\ref{eq:2.6}) 
for the PLM apparently manifest the projective duality. It results in the 
duality $f \leftrightarrow \nu$ of all formulae derived from
(\ref{eq:2.5}) and (\ref{eq:2.6}). 
\begin{lemma}
For the PLM (\ref{eq:2.5}) the relations hold
\begin{eqnarray}
\langle f_x , \nu \rangle =\langle f_x , \nu_x \rangle =
\langle f_{xx} , \nu \rangle=\langle f , \nu_{xx} \rangle =
\langle f_{xx} , \nu_{xx} \rangle=0 
\label{eq:2.7} \; ,\\
\langle f_y , \nu \rangle =\langle f_y , \nu_y \rangle =
\langle f_{yy} , \nu \rangle=\langle f , \nu_{yy} \rangle =
\langle f_{yy} , \nu_{yy} \rangle=0 \; .
\label{eq:2.8}
\end{eqnarray}
\end{lemma} 
To prove (\ref{eq:2.7}) and (\ref{eq:2.8}) we present (\ref{eq:2.5})
and their differential consequences (\ref{eq:2.9})
\beq
f \wedge f_{xx} = \star \left( \nu \wedge \nu_{xx} \right) \;\; , \;\;
f \wedge f_{yy} = -\star \left( \nu \wedge \nu_{yy} \right)
\label{eq:2.9}
\ee
in a component form:
\begin{eqnarray}
f_i f_{kx}-f_k f_{ix} & = & \eps_{iklm} \nu_l \nu_{mx}
\label{eq:2.10} \; ,\\
f_i f_{ky}-f_k f_{iy} & = & -\eps_{iklm} \nu_l \nu_{my} \;\;\;\;\;
(i,k=1,\dots,4) 
\label{eq:2.11}
\end{eqnarray}
and
\begin{eqnarray}
f_i f_{kxx}-f_k f_{ixx} & = & \eps_{iklm} \nu_l \nu_{mxx}
\label{eq:2.12} \; ,\\
f_i f_{kyy}-f_k f_{iyy} & = & -\eps_{iklm} \nu_l \nu_{myy} \;\;\;\;\;
(i,k=1,\dots,4) \;\; .
\label{eq:2.13}
\end{eqnarray}
The pairing (\ref{eq:2.1}) is obviously compatible with (\ref{eq:2.5})
and (\ref{eq:2.12}), (\ref{eq:2.13}). Equations (\ref{eq:2.10}), 
(\ref{eq:2.11}) imply 
\beq
\begin{array}{l}
f \langle f_x , \nu \rangle -
\langle f , \nu \rangle f_x = 0 \; , \\
f \langle f_y , \nu \rangle -
\langle f , \nu \rangle f_y = 0 \; , \\
f \langle f_x , \nu_x \rangle -
\langle f , \nu_x \rangle f_x = 0 \; , \\
f \langle f_y , \nu_y \rangle -
\langle f , \nu_y \rangle f_y = 0  
\end{array}
\label{eq:2.14}
\ee
while (\ref{eq:2.12}), (\ref{eq:2.13}) give 
\beq
\begin{array}{l}
f \langle f_{xx} , \nu_{xx} \rangle 
- f_{xx} \langle f , \nu_{xx} \rangle = 0 \; , \\
f \langle f_{yy} , \nu_{yy} \rangle - f_{yy} \langle f , 
\nu_{yy} \rangle = 0 \; .
\end{array}
\label{eq:2.15}
\ee
For generic $f$ the relations (\ref{eq:2.14}), (\ref{eq:2.15})
are equivalent to the relations (\ref{eq:2.7}) and (\ref{eq:2.8}).
Further from (\ref{eq:2.10}), (\ref{eq:2.11}) one gets 
\beq
f \langle f_x , \nu_y \rangle = 
f \langle f_y , \nu_x \rangle=
\left[ \nu, \nu_x, \nu_y \right] \;\; .
\label{eq:2.16}
\ee
Pairing of both sides of (\ref{eq:2.16}) with $\nu_{xy}$ and use of
(\ref{eq:2.3}) give
\beq
\langle f_x , \nu_y \rangle \langle f , \nu_{xy} \rangle=
-det|\nu, \nu_x, \nu_y, \nu_{xy}| \;\; .
\label{eq:2.17}
\ee
Since $\langle f , \nu_{xy} \rangle=
-\langle f_x , \nu_y \rangle$ one gets
\beq
\langle f_x , \nu_y \rangle^2=det|\nu, \nu_x, \nu_y, \nu_{xy}|
\;\; .
\label{eq:2.18}
\ee
Thus, using (\ref{eq:2.16}) and (\ref{eq:2.18}), we prove the
\begin{theorem}
For the PLM (\ref{eq:2.5}) $L:\nu \rightarrow f$ one has
\beq
f=\frac{\left[\nu, \nu_x, \nu_y \right]}{\sqrt{det|\nu,
\nu_x, \nu_y, \nu_{xy}|}} \;\; .
\label{eq:2.19}
\ee
For the inverse PLM $f \rightarrow \nu$ one has
\beq
\nu=\frac{\left[f, f_x, f_y \right]}{\sqrt{det| f,
f_x, f_y, f_{xy}|}} \;\; .
\label{eq:2.20}
\ee
\end{theorem}
Using now (\ref{eq:2.12}), one gets
\beq
\langle f_{xx} , \nu_x \rangle f=-\left[\nu, \nu_x, \nu_{xx} \right]
\;\; .
\label{eq:2.21}
\ee
The equality (\ref{eq:2.12}) implies that
\beq
\langle f_{xx} , \nu_x \rangle \langle f , \nu_{xxx} \rangle =
det|\nu,\nu_x, \nu_{xx}, \nu_{xxx}| \;\; .
\label{eq:2.22}
\ee
Since $\langle f , \nu_{xxx} \rangle=
 \langle f_{xx} , \nu_x $ one obtains
\beq
\langle f_{xx} , \nu_x \rangle^2=det|\nu, \nu_x, \nu_{xx}, \nu_{xxx}|
\;\; .
\label{eq:2.23}
\ee
So
\beq
f=-\frac{\left[\nu, \nu_x, \nu_{xx} \right]}{\sqrt{det|\nu,
\nu_x, \nu_{xx}, \nu_{xxx}|}} \;\; .
\label{eq:2.24}
\ee
Analogously from (\ref{eq:2.13}), one gets
\beq
\langle f_{yy} , \nu_y \rangle^2=-det|\nu, \nu_y, \nu_{yy}, \nu_{yyy}|
\label{eq:2.25}
\ee
and
\beq
f=-\frac{\left[\nu, \nu_y, \nu_{yy} \right]}{\sqrt{det|\nu,
\nu_y, \nu_{yy}, \nu_{yyy}|}} \;\; .
\label{eq:2.26}
\ee
Using the inverse PLM (\ref{eq:2.6}) (which coincides with the direct one)
, one gets
\beq
\begin{array}{lll}
\langle f_y , \nu_x \rangle^2 & = & det|\nu, \nu_x, \nu_{y}, \nu_{xy}|
\; , \\
\langle f_x , \nu_{xx} \rangle^2 & = & 
det|\nu, \nu_x, \nu_{xx}, \nu_{xxx}| \; , \\
\langle f_y , \nu_{yy} \rangle^2 & = &
-det|\nu, \nu_y, \nu_{yy}, \nu_{yyy}| \;\; .
\end{array}
\label{eq:2.27}
\ee
Comparing (\ref{eq:2.18}), (\ref{eq:2.23}), (\ref{eq:2.25}) with 
(\ref{eq:2.27}) and taking into account that
$\langle f_y , \nu_x \rangle = 
\langle f_x , \nu_y \rangle$,
$\langle f_{xx} ,
\nu_x \rangle = -\langle f_x , \nu_{xx} \rangle$, 
$\langle f_{yy} ,
\nu_y \rangle = -\langle f_y , \nu_{yy} \rangle$
one gets
\begin{theorem}
The full determinants are invariant under the PLM (\ref{eq:2.5}):
\beq
\begin{array}{l}
det|f, f_x, f_{y}, f_{xy}| = det|\nu, \nu_x, \nu_{y}, \nu_{xy}| \; , \\
det|f, f_x, f_{xx}, f_{xxx}| = det|\nu, \nu_x, \nu_{xx}, \nu_{xxx}| \; ,\\
det|f, f_y, f_{yy}, f_{yyy}| = det|\nu, \nu_y, \nu_{yy}, \nu_{yyy}| \;\; . \\
\end{array}
\label{eq:2.28}
\ee
\end{theorem}

The formulae (\ref{eq:2.27}) and (\ref{eq:2.25}) provide us the following 
expressions for projective Fubini forms (which are invariant under 
unimodular projective transformations):
\beq
\begin{array}{l}
F_2=2 \langle f_x , \nu_y \rangle dx \, dy = 2
\sqrt{det|f, f_x, f_y, f_{xy}|} dx dy =
2\sqrt{det|\nu, \nu_x, \nu_y, \nu_{xy}|} dx dy \;, 
\\
F_3= \langle f_x , \nu_{xx} \rangle dx^3 = 
\sqrt{ det|f, f_x, f_{xx}, f_{xxx}|} dx^3 =
\sqrt{ det|\nu, \nu_x, \nu_{xx}, \nu_{xxx}|} dx^3 \;, 
\\
\tilde{F}_3= \langle f_y , \nu_{y} \rangle dy^3 = 
\sqrt{ -det|f, f_y, f_{yy}, f_{yyy}|} dy^3 =
\sqrt{ -det|\nu, \nu_y, \nu_{yy}, \nu_{yyy}|} dy^3 
\;\; . 
\end{array}
\label{eq:2.29}
\ee
Further comparing (\ref{eq:2.19}), (\ref{eq:2.24}), (\ref{eq:2.26}) and 
their dual analogs, one arrives at the following
\begin{theorem}
The compatibility conditions for the PLM (\ref{eq:2.5}) are the following
\beq
\begin{array}{l}
\nu_{xx}=U_1 \nu_x + V_1 \nu_y +W_1 \nu \; , \\
\nu_{yy}=U_2 \nu_x + V_2 \nu_y +W_2 \nu 
\label{eq:2.30}
\end{array}
\ee
and
\beq
\begin{array}{l}
f_{xx}=U_1 f_x - V_1 f_y + \tilde{W}_1 f \; , \\
f_{yy}=-U_2 f_x + V_2 f_y +\tilde{W}_2 f
\label{eq:2.31}
\end{array}
\ee
where
\beq
V_1^2 = 
\frac{det|\nu, \nu_x, \nu_{xx}, \nu_{xxx}|}{det|\nu, \nu_x, \nu_{y}, 
\nu_{xy}|} \;\; , \;\;
U_2^2 = -
\frac{det|\nu, \nu_y, \nu_{yy}, \nu_{yyy}|}{det|\nu, \nu_x, \nu_{y}, 
\nu_{xy}|}
\label{eq:2.32}
\ee
and $U_1$, $W_1$, $V_2$, $W_2$, $\tilde{W}_1$, $\tilde{W}_2$ are some 
functions. 
\end{theorem}
Note that neither $\nu$ nor $f$ obey an equation of the form
$$ f_{xy}=Cf_x+Df_y+Ef \;\; . $$
Note also that $V_1$ and $U_2$ are projective invariants and have the same 
form in terms of $f$.

Equations (2.30) and (2.31) are known one. They define surfaces in the three-
dimensional projective spaces  dual to each other (see {\em e.g.} \bib{r7}).
The relations of the form (2.19),(2.20) and (2.28) derived in a different
situation also can be found  in \bib{r7}.

\section{The PLM for elliptic surfaces}
\setcounter{equation}{0}

Similar to the standard affine Lelieuvre formula for elliptic surfaces (
see {\em e.g.} \bib{r4}) there is an elliptic version of the PLM 
(\ref{eq:2.5}). 
\begin{definition}
An elliptic version of the PLM is given by the relation 
\beq
f \wedge df = \star \left( \nu \wedge \star d\nu \right)
\label{eq:3.1}
\ee
where $f \subset P^3$, 
\beq
\langle f , \nu \rangle = 0
\label{eq:3.2}
\ee
and $\star d\nu$ is a dual $1-$form. 
\end{definition}
In local coordinates $(x,y)$, $f=f(x,y)$, $\nu=\nu(x,y)$
and (\ref{eq:3.1}) is
\beq
f \wedge f_x = - \star \left( \nu \wedge \nu_y \right) \;\;, \;\;
f \wedge f_y = \star \left( \nu \wedge \nu_x \right) \;\; .
\label{eq:3.3}
\ee
The inverse PLM $L^{-1}:f \rightarrow \nu$ is given by 
$\nu \wedge \star d\nu = \star \left( f \wedge df \right) $. 
The differential consequences of (\ref{eq:3.3}) are of the form
\beq
\begin{array}{l}
f \wedge f_{xx} = - \star \left( \nu_x \wedge \nu_y  \right)
-\star \left( \nu \wedge \nu_{xy}  \right) \;\;, \\
f \wedge f_{yy} = \star \left( \nu_y \wedge \nu_x  \right)
+ \star \left( \nu \wedge \nu_{xy}  \right) 
\end{array}
\label{eq:3.4}
\ee
and
\beq
\begin{array}{l}
f_y \wedge f_{x}+f \wedge f_{xy} = - \star 
\left( \nu \wedge \nu_{yy} \right) \;, \\
f_x \wedge f_{y} + f \wedge f_{xy} = \star 
\left( \nu \wedge \nu_{xx} \right)
\;\; .\\
\end{array}
\label{eq:3.5}
\ee
Using (\ref{eq:3.2}), (\ref{eq:3.3}), one gets
\begin{lemma}
For the PLM (\ref{eq:3.1}) one has
\beq
\begin{array}{l}
\langle f , \nu_x \rangle =\langle f_x , \nu \rangle =
\langle f , \nu_y \rangle =\langle f_y , \nu \rangle = 0 \;, \\
\langle f_x , \nu_y \rangle =\langle f_y , \nu_x \rangle=0 \; ,\\
\langle f_{xy} , \nu \rangle =\langle f , \nu_{xy} \rangle=0
\end{array}
\label{eq:3.6}
\ee
and
\beq
\langle f_x , \nu_x \rangle =\langle f_y , \nu_y \rangle \;\; .
\label{eq:3.7}
\ee
\end{lemma}
Equations (\ref{eq:3.3}) imply that
\beq
f \langle f_x , \nu_x \rangle = \left[ \nu, \nu_x, \nu_y \right]
\label{eq:3.8}
\ee
and
\beq
f \langle f_y , \nu_y \rangle = \left[ \nu, \nu_x, \nu_y \right]
\label{eq:3.9}
\ee
From (\ref{eq:3.8}) one gets
\beq
\langle f , \nu_{xx} \rangle \langle f_x , \nu_x \rangle =
-det|\nu, \nu_x, \nu_{xy}, \nu_{xx}| \;\; .
\label{eq:3.10}
\ee
Since $\langle f , \nu_{xx} \rangle=-\langle f_x 
, \nu_x \rangle$,
one obtains
\beq
\langle f_x , \nu_x \rangle^2=det|\nu, \nu_x, \nu_{y}, \nu_{xx}|
\;\; .
\label{eq:3.11}
\ee
Analogously (\ref{eq:3.9}) gives
\beq
\langle f_y , \nu_y \rangle^2=det|\nu, \nu_x, \nu_{y}, \nu_{yy}|
\;\; .
\label{eq:3.12}
\ee
Thus as a consequence of (\ref{eq:3.8}), (\ref{eq:3.9}), 
(\ref{eq:3.11}), (\ref{eq:3.12}) one has
\begin{theorem}
For the PLM map (\ref{eq:3.1})
\beq
f=\frac{\left[ \nu, \nu_x, \nu_y \right]}{\sqrt{det|\nu,
\nu_x, \nu_{y}, \nu_{xx}|}}=\frac{\left[ \nu, \nu_x, \nu_y 
\right]}{\sqrt{det|\nu,
\nu_x, \nu_{y}, \nu_{yy}|}} \;\; .
\label{eq:3.13}
\ee
\end{theorem}
In a similar manner one can show that for the inverse PLM 
\beq
\nu=\frac{\left[ f, f_x, f_y \right]}{\sqrt{det|f,
f_x, f_{y}, f_{xx}|}}=\frac{\left[ f, f_x, f_y 
\right]}{\sqrt{det|f,
f_x, f_{y}, f_{yy}|}}
\label{eq:3.14}
\ee
and
\beq
\begin{array}{l}
\langle f_x , \nu_x \rangle^2= -det|f, f_x, f_{y}, f_{xx}| \;, \\
\langle f_y , \nu_y \rangle^2=-det|f, f_x, f_{y}, f_{yy}| \;\; .
\end{array}
\label{eq:3.15}
\ee
Comparison of (\ref{eq:3.11}), (\ref{eq:3.12}) with (\ref{eq:3.15})
leads to
\begin{theorem}
Full determinants change signs under the PLM (\ref{eq:3.1}):
\beq
\begin{array}{l}
det|f, f_x, f_{y}, f_{xx}| =-det|\nu, \nu_x, \nu_{y}, \nu_{xx}| \;, \\
det|f, f_x, f_{y}, f_{yy}| =-det|\nu, \nu_x, \nu_{y}, \nu_{yy}| \;\; .
\end{array}
\label{eq:3.16}
\ee
\end{theorem}
Further from (\ref{eq:3.4}) one gets 
\beq
f \langle f_{xx} , \nu_{x} \rangle 
=\left[\nu, \nu_x, \nu_{xy} \right] \;\; .
\label{eq:3.17}
\ee
Since $\langle f , \nu_{y} \rangle=0$ (\ref{eq:3.17}) implies
\beq
det|\nu, \nu_x, \nu_{y}, \nu_{xy}|=0 \;\; .
\label{eq:3.18}
\ee
Then equations (\ref{eq:3.16}) and (\ref{eq:3.11}), (\ref{eq:3.12})
($\langle f_{x} , \nu_{x} \rangle=
\langle f_{y} , \nu_{y} \rangle=0$) 
imply
\beq
det|\nu, \nu_x, \nu_{y}, \nu_{yy}|=-det|\nu, \nu_x, \nu_{y}, \nu_{xx}|
\;\; .
\label{eq:3.19}
\ee
Using (\ref{eq:3.16}), (\ref{eq:3.17}) and relations
\beq
\begin{array}{l}
f \wedge \left( f_{yy}-f_{xx} \right) = 2 \star \left( \nu
\wedge \nu_{xy} \right)
\; , \\
2 f \wedge f_{xy}= -\star \left( \nu \wedge \left[ 
\nu_{yy}-\nu_{xx} \right] \right)
\end{array}
\label{eq:3.20}
\ee
one gets
\begin{theorem}
The compatibility conditions for the PLM (\ref{eq:3.1}) are of the
form
\begin{eqnarray}
\nu_{xy} & = & U \nu_x + V \nu_y +W \nu \;,  
\label{eq:3.21} \\
\nu_{yy} & - & \nu_{xx} = -2 \tilde{V} \nu_x + 2 \tilde{U} \nu_y +C \nu 
\label{eq:3.22} 
\end{eqnarray}
and
\begin{eqnarray}
f_{xy} & = & \tilde{U} f_x + \tilde{V} f_y + \tilde{W} f \;,  
\label{eq:3.23} \\
f_{yy} & - & f_{xx} = -2 V f_x + 2 U f_y +\tilde{C} f 
\label{eq:3.24} 
\end{eqnarray}
where $U$, $V$, $W$, $C$, $\tilde{U}$, $\tilde{V}$, $\tilde{W}$,
$\tilde{C}$ are some functions.
\end{theorem}

Equations (\ref{eq:3.23}), (\ref{eq:3.24}) characterize a surface in $P^3$ 
parameterized by conjugate lines (see \bib{r7}). 

Thus, the PLM (\ref{eq:3.1}) is the map between surfaces in dual spaces 
$P_3$ and $P^3$ parameterized by conjugate lines. 

The PLM's (\ref{eq:2.5}) and (\ref{eq:3.1}) and the corresponding formulae 
can be written in a unique common form. For this purpose we introduce the 
variables $\xi$ and $\eta$ defined as $\xi=x$, $\eta=y$ in the case
(\ref{eq:2.5}) and as $\xi=x+iy=z$, $\eta=x-iy=\zba$ in the case
(\ref{eq:3.1}). Then the formulae (\ref{eq:2.5}) and (\ref{eq:3.1})
take the form
\beq
f \wedge f_\xi = \sigma \star \left( \nu \wedge \nu_\xi \right) \;\; ,
\;\;
f \wedge f_\eta = -\sigma \star \left( \nu \wedge \nu_\eta \right)
\label{eq:3.25}
\ee
where $\sigma=1$ in the real case and $\sigma=-\sqrt{-1}$, $\xi=z$,
$\eta=\zba$ for the case considered in this section. 

Note that there are 
other compatibility conditions for the formulae (\ref{eq:3.25}) different 
from those given by (\ref{eq:2.29}), (\ref{eq:2.30}) or 
(\ref{eq:3.21})-(\ref{eq:3.24}). Indeed, written in coordinates formulae 
(\ref{eq:3.25}) are equivalent to the following
\beq
\left( \frac{f_k}{f_i} \right)_\xi = \sigma f_i^{-2} \eps_{iklm} \nu_l 
\nu_{m\xi} \; \; , \;\;
\left( \frac{f_k}{f_i} \right)_\eta =- \sigma f_i^{-2} \eps_{iklm} \nu_l 
\nu_{m\eta} \;\; . 
\label{eq:3.26}
\ee
An obvious compatibility condition for (\ref{eq:3.26}) is equivalent to the 
system 
\beq
\nu_{k\xi\eta}=\left( \log{f_i} \right)_\eta \nu_{k\xi}+
\left( \log{f_i} \right)_\xi \nu_{k\eta} +u_i \nu_k
\label{eq:3.27}
\ee
($i \neq k$, $i, k= 1,\ldots,4$) where $u_i$ are some function. 

From the inverse formulae (\ref{eq:3.25}) one gets
\beq
f_{k\xi\eta}=\left( \log{\nu_i} \right)_\eta f_{k\xi}+
\left( \log{\nu_i} \right)_\xi f_{k\eta} +\tilde{u}_i f_k
\;\;\;\; (i \neq k)
\label{eq:3.28}
\ee
where $\tilde{u}_i$ are some function. 
However, equations (\ref{eq:3.27}) and (\ref{eq:3.28}) are not
form-invariant under projective transformations. So in contrast to 
equations (\ref{eq:2.29}), (\ref{eq:2.30}) or
(\ref{eq:3.21})-(\ref{eq:3.24}), they do not characterize projective 
properties of surfaces. 

All the results derived until this 
section for the real projective space $RP^3$ are apparently extendable
to the space $CP^3$. In this case $f \subset CP^3$, $\nu \subset
CP_3$ and $\xi$, $\eta$ are complex 
valued variables. An intermediate case $f \subset CP^3$, 
$\nu \subset CP_3$ and $\xi=z$, 
$\eta=\zba$, $z \in C$, $\sigma=-\sqrt{-1}$ which provides the PLM for 
surfaces with positive-defined metric in $CP^3$ could be of particular 
interest to the theory of Riemann surfaces.

The PLM discussed above was formulated in asymptotic coordinates or in 
conjugate line coordinates. The PLM for surfaces in $RP^3$ can be 
formulated in general coordinates on the surfaces. We will get these 
formulae in the next section as a particular case of the PLM for 
hypersurfaces.

\section{The projective Lelieuvre map for hypersurfaces.}
\setcounter{equation}{0}

Let $f=(f_1, \ldots, f_{n+2})$ and $\nu=(\nu_1, \ldots, \nu_{n+2})$ 
be homogeneous coordinates in dual projective spaces 
$P^{n+1}$ and $P_{n+1}$ paired by (\ref{eq:2.1}). Consider hypersurfaces 
$M:f(x_1,\ldots,x_n) \subset P^{n+1}$ and 
$M^*:\nu(x_1,\ldots,x_n) \subset P_{n+1}$
where $x_1,\ldots,x_n$ are any local coordinates on surfaces.
\begin{definition}
The PLM $L:\nu \rightarrow f$
for hypersurfaces in $P^{n+1}$ is defined by the system of 
equation
\beq
f \wedge f_{x_{\alpha}}= \sum_{\beta =1}^{n}
A_{\alpha \beta} \star \left( 
\nu_{x_{1}} \wedge \ldots \wedge \nu_{x_{\beta-1}} \wedge \nu
\wedge \nu_{x_{\beta+1}} \wedge \ldots \wedge \nu_{x_{n}} \right) \;\;,
\;\; \alpha=1,\ldots,n
\label{eq:4.1}
\ee
where $A_{\alpha \beta}$ ($\alpha$, $\beta=1,\ldots,n$) are functions on
$x_1,\ldots,x_n$.
\end{definition}
The inverse PLM $L^{-1}:f \rightarrow \nu$ is given by 
\beq
\nu_{x_{1}} \wedge \ldots \wedge \nu_{x_{\beta-1}} \wedge \nu
\wedge \nu_{x_{\beta+1}} \wedge \ldots \wedge \nu_{x_{n}}=
\star \sum_{\gamma=1}^{n} \left( A^{-1} \right)_{\beta \gamma} f
 \wedge f_{x_{\gamma}} \;\;, \;\; \beta=1,\ldots,n
\label{eq:4.2}
\ee
where $A^{-1}$ is the matrix inverse to the matrix $A$ ($\det A \neq
0$). 

In local coordinates in $P_{n+1}$ the formulae (\ref{eq:3.1}) look like
\beq
\begin{array}{l}
f_i f_{kx_\alpha}-f_k f_{ix_\alpha}=
\sum_{\beta =1}^{n} A_{\alpha \beta}\eps_{ikl_1 \ldots l_n}
\nu_{l_1 x_1} \dots \nu_{l_{\beta-1} x_{\beta-1}}
\nu_{l_{\beta}} \nu_{l_{\beta+1} x_{\beta+1}} \ldots 
\nu_{l_{n} x_{n}} \; , \\
(i,k=1,\dots,n+2\;\;,\;\; \alpha=1,\ldots,n) \;\;.
\end{array}
\label{eq:4.3}
\ee
\begin{lemma}
For the PLM (\ref{eq:4.1}) one has 
\beq
\langle f_{x_{\alpha}} , \nu \rangle 
=\langle f , \nu_{x_{\alpha}} \rangle=0 
\;\;,\;\; \alpha=1,\ldots,n \;\;.
\label{eq:4.4}
\ee
\end{lemma}
Equations (\ref{eq:4.3}) imply that
\beq
\langle f , \nu \rangle f_{x_{\alpha}}-\langle f_{x_{\alpha}}
, \nu
\rangle f=0\;\;,\;\; \alpha=1,\ldots,n \;\; .
\label{eq:4.5}
\ee
Due to (\ref{eq:2.1}), one gets (\ref{eq:4.4}).

Further (\ref{eq:4.1}) implies
\beq
\langle f_{x_\alpha} , \nu_{x_\gamma} \rangle f =
\sum_{\beta=1}^{n} A_{\alpha\beta} \left[
\nu_{x_\gamma}, \nu_{x_1}, \ldots, \nu_{x_{\beta-1}},\nu,
\nu_{x_{\beta+1}}, \ldots, \nu_{x_n} \right]
\;\;, \;\; \alpha, \gamma=1,\ldots,n \;\; .
\label{eq:4.6}
\ee
Since
\beq
\left[
\nu_{x_\gamma}, \nu_{x_1}, \ldots, \nu_{x_{\beta-1}},\nu,
\nu_{x_{\beta+1}}, \ldots, \nu_{x_n} \right] =
-\delta_{\gamma \beta}
\left[\nu, \nu_{x_1},\ldots, \nu_{x_n} \right]
\label{eq:4.7}
\ee
where $\delta_{\alpha \beta}$ is the Kroneker symbol, one has
\beq
\langle f_{x_\alpha} , 
\nu_{x_\gamma} \rangle f =-A_{\alpha\gamma}
\left[\nu, \nu_{x_1},\ldots, \nu_{x_n} \right]
\;\;,\;\; \alpha,\gamma=1,\ldots,n\;\;.
\label{eq:4.8}
\ee
It follows from (\ref{eq:4.8}) that 
\beq
\langle f_{x_\alpha} , \nu_{x_\gamma} \rangle
\langle f_{x_\beta} , \nu_{x_\delta} \rangle =
-A_{\alpha\gamma} det|\nu_{x_\beta x_\delta}, \nu,
\nu_{x_1}, \ldots,\nu_{x_n}| \;\;,\;\; \alpha,\beta,\delta,\gamma=1,\ldots,
n \;.
\label{eq:4.9}
\ee
Exchanging indices $(\alpha,\beta) \leftrightarrow (\beta, \alpha)$ in
(\ref{eq:4.9}), one also gets
\beq
\langle f_{x_\beta} , \nu_{x_\delta} \rangle
\langle f_{x_\alpha} , \nu_{x_\gamma} \rangle =
-A_{\beta\delta} det|\nu_{x_\alpha x_\gamma}, \nu,
\nu_{x_1}, \ldots,\nu_{x_n}| \;\;,\;\; \alpha,\beta,\delta,\gamma=1,\ldots,
n \;.
\label{eq:4.10}
\ee
From (\ref{eq:4.9}) or (\ref{eq:4.10}) it follows at $\alpha=\beta$,
$\gamma=\delta$ that
\beq
\langle f_{x_\alpha} , \nu_{x_\gamma} \rangle^2=
A_{\alpha\gamma} det|\nu_{x_\alpha x_\gamma}, \nu, 
\nu_{x_1}, \ldots,\nu_{x_n}| \;\;.
\label{eq:4.11}
\ee
Taking into account (\ref{eq:4.8}) and (\ref{eq:4.11}), one gets the 
following
\begin{theorem}
For the PLM (\ref{eq:4.1}) for hypersurfaces one has 
\beq
f=- \left( \frac{A_{\alpha\gamma}}{det|\nu_{x_\alpha x_\gamma},
\nu,\nu_{x_1},\ldots,\nu_{x_n}|} \right)^{\frac{1}{2}}
\left[ \nu,\nu_{x_1},\ldots,\nu_{x_n} \right] \;\; .
\label{eq:4.12}
\ee
\end{theorem}
Further comparing (\ref{eq:4.9}) and (\ref{eq:4.10}), one obtains the 
equation
\beq
det|A_{\alpha\gamma} \nu_{x_\beta x_\delta}-
A_{\beta\delta} \nu_{x_\alpha x_\gamma}, \nu,\nu_{x_1},\ldots,\nu_{x_n}|
=0 \;\;,\;\; \alpha,\beta,\delta,\gamma=1,\ldots,n\;\;.
\label{eq:4.13}
\ee
This equation implies the following
\begin{theorem}
The compatibility conditions for the PLM (\ref{eq:4.1}) are given by the 
system of equations
\beq
A_{\alpha\gamma} \nu_{x_\beta x_\delta}-
A_{\beta\delta} \nu_{x_\alpha x_\gamma}+
\sum_{\rho=1}^{n} {U_{\alpha,\beta,\gamma,\delta}^{(\rho)} 
\nu_{x_\rho}} +W_{\alpha,\beta,\gamma,\delta} \nu=0
\;\;,\;\;\alpha,\beta,\gamma,\delta=1,\ldots,n
\label{eq:4.14}
\ee
$U_{\alpha,\beta,\gamma,\delta}^{(\rho)}$
and $W_{\alpha \beta \gamma \delta}$ are some functions. 
\end{theorem}
These functions vanish when simultaneously $\alpha=\beta$ and
$\gamma=\delta$ and for those ${\alpha, \beta, \gamma, \delta}$
for which both $A_{\alpha \gamma}=0$ and $A_{\beta \delta}=0$.
\begin{corollary}
In virtue of the compatibility condition (\ref{eq:4.14}) the factor in the 
formula (\ref{eq:4.12}) is independent on choice of indices $\alpha$,
$\beta$.
\end{corollary}

In the particular case $n=2$ the formula (\ref{eq:4.1}), (\ref{eq:4.12})
(\ref{eq:4.14}) give the PLM for surface in $P_3$ ($P^3$) in general 
coordinate system. At $A_{11}=A_{22}=0$, $A_{12}=A_{21}=-2$ one reproduces
the results of section $3$ while at the case $A_{12}=A_{21}=0$,
$A_{11}=A_{22}=-2$ one gets the formulae of section $3$.

\section{Projective Lelieuvre map for discrete surfaces}
\setcounter{equation}{0}

Discrete surfaces (maps $\bfZ_2 \rightarrow R^N$) 
are the subject of intensive study now (see 
{\em e.g.} \bib{r8}-\bib{r9}). A discrete analog of the Lelieuvre formula 
for discrete affine spheres has been found recently in \bib{r10}. 

Here we present the projective Lelieuvre formulae for discrete surfaces 
in $P^3$. So let $f:{\bfZ}_2 \rightarrow P^3$ and 
$\nu:{\bfZ}_2 \rightarrow P_3$ with the pairing (\ref{eq:2.1}).
Thus $f=f(n_1,n_2)$ and $\nu=\nu(n_1,n_2)$
where $n_1$, $n_2$ are integers. We denote the shift 
operators as $T_1$ and $T_2$: $T_1 f(n_1,n_2)=f(n_1+1,n_2)$, 
$T_2 f(n_1,n_2)=f(n_1,n_2+1)$. 
For compactness we will denote $f_1=T_1 f$, $f_{11}=T_{1}^2 f$,
$f_2=T_2 f$, $f_{22}=T_2^2 f$, $f_{-1}=T^{-1}_1 f$ etc. and will omit 
arguments of $f$ and $\nu$.
\begin{definition}
Discrete PLM $L:\nu \rightarrow f$ is given by relations
\beq
f \wedge f_1 = \star \left( \nu \wedge \nu_1 \right) \;\; , \;\;
f \wedge f_2 = - \star \left( \nu \wedge \nu_2 \right) \;\; .
\label{eq:5.1}
\ee
\end{definition}
The inverse map $f \rightarrow \nu$ is of the same form
\beq
\nu \wedge \nu_1= \star \left( f \wedge f_1 \right) \;\; , \;\;
\nu \wedge \nu_2=-\star \left( f\wedge f_2\right) \;\; .
\ee
In coordinates (\ref{eq:5.1}) looks like
\beq
\begin{array}{l}
f_i f_{1k}-f_k f_{1i}=\frac{1}{2} \eps_{iklm} \left( \nu_l \nu_{1m}
-\nu_m \nu_{1l} \right) \; , \\
f_i f_{2k}-f_k f_{2i}=-\frac{1}{2}\eps_{iklm} \left( \nu_l\nu_{2m}
-\nu_m \nu_{2l} \right) 
\end{array}
\;\;\;\; i,k=1,\ldots,4 \;\; .
\label{eq:5.2}
\ee
So the PLM (\ref{eq:5.1}) is, in fact, the identification:
{\em 1)} of the polar Pl\"{u}cker coordinates of discrete surface in
$P_3$ in direction $T_1^n \nu$
with the corresponding Pl\"{u}cker coordinates in $P^3$ and
{\em 2)} of the anti-polar Pl\"{u}cker coordinates in $P_3$ in direction 
$T_2^n \nu$ with the corresponding Pl\"{u}cker coordinates in $P^3$.
\begin{lemma}
For discrete PLM (\ref{eq:5.1}) one has
\beq
\langle f_{\alpha} , \nu \rangle =
\langle f , \nu_{\alpha} \rangle =0
\;\;, \;\; \alpha=1,2 \;\; .
\label{eq:5.3}
\ee
\end{lemma}
From (\ref{eq:5.2}) one gets 
\beq
\begin{array}{l}
f \langle f_{\alpha} , \nu \rangle - 
\langle f , \nu \rangle f_{\alpha} =0 \;
, \\
\langle f_{\alpha} , \nu_{\alpha}\rangle- 
\langle f , \nu_{\alpha}\rangle
f_{\alpha} =0 \;\; , \;\; \alpha=1,2 \;\; .
\end{array}
\label{eq:5.4}
\ee
Since $\langle f , \nu \rangle=
\langle f_{\alpha} , \nu_{\alpha}\rangle = 0$ 
one gets (\ref{eq:5.3}).

The relations (\ref{eq:5.2}) and their shifted versions give
\beq
f \langle f_1 , \nu_2 \rangle 
= \left[ \nu,\nu_1,\nu_2 \right] \;\;, 
\;\;f \langle f_2 , \nu_1 \rangle =
\left[ \nu,\nu_1, \nu_2 \right] \;, 
\label{eq:5.5}
\ee
\beq
f_1 \langle f_{11} , \nu \rangle 
= \left[ \nu,\nu_1,\nu_{11} \right] \;\;, 
\;\; f_1 \langle f_{12} , \nu \rangle =
\left[ \nu,\nu_1, \nu_{12} \right] \;, 
\label{eq:5.6}
\ee
\beq
f_2 \langle f_{22} , \nu \rangle =
- \left[ \nu,\nu_2,\nu_{22} \right] \;\;, 
\;\; f_2 \langle f_{12} , \nu \rangle =
\left[ \nu,\nu_2, \nu_{12} \right] \;, 
\label{eq:5.7}
\ee
and
\beq
f_{12} \langle f_{2} , \nu_1 \rangle 
=- \left[ \nu,\nu_2,\nu_{12} \right] \;\; .
\label{eq:5.8}
\ee
From (\ref{eq:5.5})-(\ref{eq:5.8}) it follows
\begin{lemma}
For the discrete PLM (\ref{eq:5.1}) one has
\beq
\langle f_{1} , \nu_2 \rangle =
\langle f_{2} ,  \nu_1 \rangle \; ,
\label{eq:5.9}
\ee
\beq
\langle f , \nu_{12} \rangle =\langle f_{12} , \nu \rangle
\label{eq:5.10}
\ee
and
\beq
\begin{array}{l}
\langle f_{12} , \nu \rangle \langle f_{1} , \nu_2 \rangle =
-det|\nu, \nu_1, \nu_{2}, \nu_{12}| \; , \\
\langle f_{11} , \nu \rangle \langle f_{1} , \nu_2 \rangle =
det|\nu, \nu_1, \nu_{2}, \nu_{11}| \; , \\
\langle f_{22} , \nu \rangle \langle f_{1} , \nu_2 \rangle =
det|\nu, \nu_1, \nu_{2}, \nu_{22}| \;\; . 
\end{array}
\label{eq:5.11}
\ee
\end{lemma}
Further equations (\ref{eq:5.5})-(\ref{eq:5.7}) give rise to the following
\begin{theorem}
The compatibility conditions for the PLM (\ref{eq:5.1}) have the form
\beq
\begin{array}{l}
\nu_{11}=A_1 \nu_{12} + B_1 \nu_1 +C_1 \nu \; , \\
\nu_{22}=A_2 \nu_{12} + B_2 \nu_2 +C_2 \nu 
\end{array}
\label{eq:5.12}
\ee
where
\beq
A_1=-\frac{\langle f_{11} , \nu \rangle}{
\langle f_{12} , \nu \rangle} 
\;\, , \;\,
C_1=\frac{\langle f_{11} , \nu_{12} 
\rangle}{\langle f_{12} , \nu \rangle} \;\; , \;\,
A_2=-\frac{\langle f_{22} , \nu 
\rangle}{\langle f_{12} , \nu \rangle} \;\; , \;\,
C_2=-\frac{\langle f_{22} , \nu_{12} 
\rangle}{\langle f_{12} , \nu \rangle} \;\;.
\label{eq:5.13}
\ee
and correspondingly
\beq
\begin{array}{l}
f_{11}=-A_1 f_{12} +\tilde{B}_1 f_1 +C_1 f \; , \\
f_{22}=-A_2 f_{12} +\tilde{B}_2 f_2 +C_2 f 
\end{array}
\label{eq:5.14}
\ee
where $B_1$, $B_2$, $\tilde{B}_1$, $\tilde{B}_2$ are some functions.
\end{theorem}
For the inverse PLM one gets the formulae (\ref{eq:5.9})-(\ref{eq:5.11})
with the substitution $f \leftrightarrow \nu$, in particular,
\beq
\nu \langle f_2 , \nu_1 \rangle= \left[ f, f_1, f_2 \right]
\ee
and
\beq
\langle f , \nu_{12} \rangle 
\langle f_2 , \nu_1 \rangle=
-det|f, f_1, f_{2}, f_{12}| \;\; .
\label{eq:5.15}
\ee
Comparing (\ref{eq:5.11}) and (\ref{eq:5.15}), one gets
\begin{theorem}
For the PLM (\ref{eq:5.1}) $\nu \rightarrow f$ one has
\beq
det|f, f_1, f_{2}, f_{12}|=det|\nu, \nu_1, \nu_{2}, \nu_{12}|
\;\; .
\label{eq:5.16}
\ee
\end{theorem}
So the volume of simplex with with vertices at the origin and the points
$\nu$, $\nu_1$, $\nu_2$, $\nu_{12}$
is preserved by the PLM (\ref{eq:5.1}).

The formulae of this section are apparently reduced to those of section 
$2$ in the continuous limit $T_{\alpha} \nu=\nu+
\frac{\partial \nu}{\partial x^{\alpha}} dx^{\alpha}$, $\alpha=1,2$,
$dx^{\alpha} \rightarrow 0$.
In particular
\beq
\begin{array}{l}
det|f, f_1, f_{2}, f_{12}| \rightarrow
det|f, f_x, f_{y}, f_{xy}| \left( dx \, dy \right)^2
\; ,  \\
det|f, f_1, f_{11}, f_{111}| \rightarrow
det|f, f_x, f_{xx}, f_{xxx}| \left( dx \right)^6
\; ,  \\
det|f, f_2, f_{22}, f_{222}| \rightarrow
det|f, f_y, f_{yy}, f_{yyy}| \left( dy \right)^6
\end{array}
\label{eq:5.17}
\ee
and similar for determinants with $\nu$.

Comparing (\ref{eq:5.17}) with (\ref{eq:2.29}) and
using (\ref{eq:5.16}), one gets
\begin{proposition}
The quantities
\beq
\begin{array}{l}
F^d_2=2 \sqrt{ \langle f_{12} , \nu \rangle 
\langle f_{1} , \nu_2 \rangle } 
=\sqrt{ det|f, f_1, f_{2}, f_{12}|}
=\sqrt{ det|\nu, \nu_1, \nu_{2}, \nu_{12}|} \; , \\
F^d_3= \sqrt{ -\langle f_{1}  \nu_{111} \rangle 
\langle f_{11} , \nu \rangle }
=\sqrt{ det|f, f_1, f_{11}, f_{111}|}
=\sqrt{ det|\nu, \nu_1, \nu_{11}, \nu_{111}|} \; , \\
\tilde{F}^d_3= \sqrt{ \langle f_{2} , \nu_{222} \rangle 
\langle f_{22} , \nu \rangle }
=\sqrt{ det|f, f_2, f_{22}, f_{222}|}
=\sqrt{ det|\nu, \nu_2, \nu_{22}, \nu_{222}|}
\end{array}
\label{eq:5.18}
\ee
are the discrete analogs of the Fubini's form (\ref{eq:2.29}).
\end{proposition}

\section{Affine and dual affine gauges}
\setcounter{equation}{0}

Let us consider an affine "reduction" of the formulae derived. To get to 
affine geometry relations one should pass to inhomogeneous coordinates
(say $\left( \frac{f_1}{f_4}, \frac{f_2}{f_4}, \frac{f_3}{f_4}\right)$)
or choose the "gauge" $f_4=-1$. There is also a possibility to do the 
same in the dual space $P_3$. 

Let us consider first the gauge $f_{4}=-1$.
In this case $\nu_4=\langle \bff , \bfnu \rangle$,
where we denote ${\bff}= \left( f_1, f_2, f_3 \right)$,
${\bfnu}= \left( \nu_1, \nu_2, \nu_3 \right)$. 
The formulae (\ref{eq:2.5}) are reduced obviously to (\ref{eq:1.1}).
Equations (\ref{eq:3.27}) with $i=4$ are affine form-invariant. So in the 
gauge $f_4=-1$ the affine conormal ${\bfnu}$ in addition to equations
(\ref{eq:2.29}) obeys also the equation
\beq
{\bfnu}_{xy}=U_4 {\bfnu} \;\; .
\label{eq:6.1}
\ee
Taking into account this equation, one can show that
\beq
det|\nu, \nu_x, \nu_{y}, \nu_{xy}| = 
\langle {\bff}_x , {\bfnu}_y \rangle
det|{\bfnu}, {\bfnu}_{x}, {\bfnu}_{y}| \;\; .
\ee
Using this relation, one obtains from (\ref{eq:2.18}) that
\beq
\langle {\bff}_x , {\bfnu}_y \rangle =
det|{\bfnu}, {\bfnu}_{x}, {\bfnu}_{y}|
\label{eq:6.2}
\ee
and hence
\beq
det|\nu, \nu_x, \nu_{y}, \nu_{xy}|=
\left( det|{\bfnu}, {\bfnu}_{x}, {\bfnu}_{y}| 
\right)^2 \;\; .
\label{eq:6.3}
\ee
Then the fourth component of the relation (\ref{eq:2.19}) is
satisfied identically while for ${\bff}$
one gets a standard expansion of ${\bff}$
in normal and tangent components. The relation (\ref{eq:2.20}) is reduced 
to the known expression of 
${\bfnu}$ in terms of ${\bff}$. The formulae
(\ref{eq:2.22}) and (\ref{eq:2.25}) give rise to
\beq
\langle {\bff}_{xx} , {\bfnu}_{x} \rangle =
det|\bfnu, \bfnu_{x}, \bfnu_{xx}| \;\;, \;\;
\langle {\bff}_{yy} , {\bfnu}_{y} \rangle =-
det|\bfnu, \bfnu_{y}, \bfnu_{yy}| \;\; .
\label{eq:6.4}
\ee
Finally the formulae (2.28) become
\beq
\begin{array}{l}
det|\bff_x,\bff_y,\bff_{xy}|=
\left( det|\bfnu,\bfnu_x,\bfnu_{xy}| \right)^2 =F^2 \;, \\
det|\bff_x,\bff_{xx},\bff_{xxx}|=
\left( det|\bfnu,\bfnu_x,\bfnu_{xx}| \right)^2 =A^2 \;, \\
det|\bff_y,\bff_{yy},\bff_{yyy}|=-
\left( det|\bfnu,\bfnu_y,\bfnu_{yy}| \right)^2 =-B^2 \;\; .
\end{array}
\label{eq:6.5}
\ee
As a result, the Fubini's forms (\ref{eq:2.29}) are reduced to
\beq
F_2 \rightarrow F \; dx \, dy \;\;, \;\;
F_3 \rightarrow A \; dx^3 \;\;, \;\;
\tilde{F}_3 \rightarrow B \; dy^3
\nonumber
\ee
{\em i.e.} to the well-known affine Blaschke and cubic Fubini's forms 
in terms of coordinates ${\bff}$
and affine conormal ${\bfnu}$ (see {\em e.g.}
\bib{r2}-\bib{r3}).

So, in the gauge $f_4=-1 $ one reproduces the relations for the affine 
Lelieuvre map for surfaces with indefinite metric.

Now let $\nu_4=1$.  So $f_4=-\langle \bff , \bfnu \rangle$.
The inverse map (\ref{eq:2.5}) $f \rightarrow \nu$
is of the form (\ref{eq:1.1}) while the direct map $\nu \rightarrow f$
is given by
\beq
\left( \frac{{\bff}}{\langle 
{\bff} , {\bfnu} \rangle } \right)_x
= \frac{ [{\bfnu}, {\bfnu}_x ]}{\left( 
\langle {\bff} , {\bfnu} \rangle  \right)^2} \;\; , \;\;
\left( \frac{{\bff}}{\langle 
{\bff} , {\bfnu}\rangle } \right)_y
=- \frac{|{\bfnu}, {\bfnu}_y]}{\left( 
\langle {\bff} {\bfnu} \rangle  \right)^2} \;\;.
\label{eq:6.6}
\ee
The compatibility condition for (\ref{eq:6.6}) are given by
(\ref{eq:2.29}), (\ref{eq:2.30}) plus equations
\beq
\begin{array}{lll}
{\bfnu}_{xy} & = & \left( \log{a} \right)_x 
{\bfnu}_y + \left( \log{a} \right)_y {\bfnu}_x
+u {\bfnu}
\; , \\
{\bff}_{xy} & = & \tilde{u} {\bff}
\end{array}
\label{eq:6.7}
\ee
where $a=-\langle {\bff} , {\bfnu} \rangle$
and $u$, $\tilde{u}$
are some functions. The vector ${\bfnu}$ is not the 
standard affine conormal but 
$\frac{{\bfnu}}{\langle {\bff} , {\bfnu} \rangle}$
is. Further the formula (\ref{eq:2.19}) gives
\beq
{\bff} = 
\frac{ [{\bfnu}_x , {\bfnu}_y]}{\sqrt{ 
det|{\bfnu}_y, {\bfnu}_x, {\bfnu}_{xy}| }}
\label{eq:6.8}
\ee
and the relations (\ref{eq:2.28}) are reduced to those of (\ref{eq:6.5})
with substitution $f \leftrightarrow \nu$, in particular, to 
\beq
\left( det|{\bff}, {\bff}_x, {\bff}_{y}|
\right)^2 = -
det|{\bfnu}_x, {\bfnu}_y, {\bfnu}_{xy}|
\label{eq:6.9}
\ee
Thus in the dual affine gauge $\nu_4=1$
the PLM generates a class of surfaces in $P_3$. 

The choice of the gauge $f_4=-1$
imposes no constraints on surfaces. The choice 
of particular gauge only destroys the projective covariance of formulae and 
symmetry between dual spaces $P_3$ and $P^3$.

If one now demands that $f_4=-1$ and $\nu_4=1$
simultaneously then one constraints surfaces by the condition
\beq
\langle {\bff} , {\bfnu} \rangle =1
\label{eq:6.10}
\ee
Such a class consists of affine spheres.
For affine spheres all formulae are symmetric under substitution
${\bff} \leftrightarrow {\bfnu}$
(see also \bib{r10}). So the affine spheres form 
the particular class of surfaces for which the general projective duality 
via the PLM (\ref{eq:2.5}) is restored as the duality on the affine level.

Similar results are valid for affine "reduction" of the PLM
(\ref{eq:3.1}).
For the PLM (\ref{eq:4.1}) for hypersurfaces the choice of the gauge
$f_{n+2}=-1$ reduces (\ref{eq:4.1}) to the formula
\beq
{\bff}_{x\alpha} =-
\sum_{\beta=1}^{n} {A_{\alpha\beta} 
\left[ {\bfnu}_{x_1},\ldots,\bfnu_{x_{\beta-1}},\bfnu,
\bfnu_{x_{\beta+1}},\ldots,\bfnu_{x_n}  \right]} \;\; .
\label{eq:6.11}
\ee
derived in \bib{r5} (${\bfnu}=(\nu_1, \ldots, \nu_{n+1})$).
From the formula (\ref{eq:4.8}) for the $n+2-$th component one gets
\beq
A_{\alpha \gamma}=-
\frac{\langle {\bff}_{x\alpha} , {\bfnu}_{x\gamma} 
\rangle}{det|{\bfnu},
{\bfnu}_{x_1},\ldots,{\bfnu}_{x_n}|} 
\label{eq:6.12}
\ee
that coincides with that of \bib{r5}. The compatibility conditions for the 
map (\ref{eq:6.11}) are given by the equations (\ref{eq:4.14}) for
${\bfnu}$.

At last, let us consider the affine gauges for the discrete PLM
(\ref{eq:5.1}). At the gauge $f_4=-1$ the formulae (\ref{eq:5.1})
take the form
\beq
{\bff}_1-{\bff}=\left[ {\bfnu},
{\bfnu}_{1} \right] \;\;, 
{\bff}_2-{\bff}=-\left[ {\bfnu},{\bfnu}_{2} \right] \;\; . 
\label{eq:6.13}
\ee
The conditions (\ref{eq:5.3}) and (\ref{eq:5.9}) become
($\nu_4=\langle \bff , \bfnu \rangle$)
\beq
\langle ( {\bff}-{\bff}_{\alpha}) ,
{\bfnu}_{\alpha} \rangle =0 \;\;, \;\;
\langle ( {\bff}_{\alpha} - {\bff}) ,
({\bfnu}_{\alpha}-{\bfnu}) \rangle =0 \;\;, \;\;
\alpha=1,2
\label{eq:6.14}
\ee
and
\beq
\langle ( {\bff}_{1}-{\bff}_{2}) ,
({\bfnu}_{1}+{\bfnu}_{2} ) \rangle =0 \;\; .
\label{eq:6.15}
\ee
which are obviously satisfied due to (\ref{eq:6.13}). 

The compatibility conditions for (\ref{eq:6.13}) are given by
equations (\ref{eq:5.12}) for ${\bfnu}$ and also by the equation
\beq
{\bfnu}_{12}+{\bfnu}=H\left({\bfnu}_1+ {\bfnu}_2 \right)
\label{eq:6.16}
\ee
where $H$ is some function. Note that equation (\ref{eq:6.16})
is not an affine reduction of some projectively covariant compatibility 
conditions. At the gauge $f_4=-1$, using (\ref{eq:5.3}) and
(\ref{eq:5.5})-(\ref{eq:5.7}), one obtains
\beq
\begin{array}{lll}
\langle f_2 , \nu_1 \rangle & = & 
\langle ( {\bff}_{1}-{\bff}) ,
({\bfnu}_{1}-{\bfnu} ) \rangle 
= det|{\bfnu},{\bfnu}_1,{\bfnu}_2| \; , \\
\langle f_{11} , \nu \rangle & = & -
\langle ( {\bff}_{11}-{\bff}) ,
({\bfnu}_{1}-{\bfnu} ) \rangle 
= det|{\bfnu},{\bfnu}_1,{\bfnu}_{11}| \; , \\
\langle f_{22} , \nu \rangle & = & -
\langle ( {\bff}_{22}-{\bff}) ,
({\bfnu}_{2}-{\bfnu} ) \rangle 
=- det|{\bfnu},{\bfnu}_2,{\bfnu}_{22}| \; , 
\end{array}
\label{eq:6.17}
\ee
and
\beq
\langle f_{22} , \nu \rangle =\langle ( {\bff}_{12}-
{\bff}) , {\bfnu} \rangle \langle( {\bff}_{12}
-{\bff}) , {\bfnu}_{12} \rangle= -
det|{\bfnu},{\bfnu}_2,{\bfnu}_{12}| \;\; .
\label{eq:6.18}
\ee
Further, using the properties of determinants and (\ref{eq:6.14}),
(\ref{eq:6.18}) at $f_4=-1$ one gets
\beq
det|f, f_1, f_{2}, f_{12}| =
det|{\bff}_{1}-{\bff}, {\bff}_{2}-{\bff},
{\bff}_{12}-{\bff} |
\label{eq:6.19}
\ee
and
\beq
det|\nu, \nu_1, \nu_{2}, \nu_{12}| = \langle ({\bff}_{1}-
{\bff}) , {\bfnu}_{12} \rangle 
det|{\bfnu},{\bfnu}_1,{\bfnu}_2|=
det|{\bfnu},{\bfnu}_1,{\bfnu}_{12}|
det|{\bfnu},{\bfnu}_1,{\bfnu}_2| \;\; .
\label{eq:6.20}
\ee
In virtue of (\ref{eq:5.16}) and (\ref{eq:6.19}), (\ref{eq:6.20})
we have 
\begin{theorem}
For the discrete affine Lelieuvre map (\ref{eq:6.13})
\beq
det|{\bff}_{1}-{\bff}, {\bff}_{2}-{\bff},{\bff}_{12}-{\bff} |=
det|{\bfnu},{\bfnu}_1,{\bfnu}_{12}|
det|{\bfnu},{\bfnu}_1,{\bfnu}_2| \;\; .
\label{eq:6.21}
\ee
\end{theorem}
In continuous limit $\bfnu_1=\bfnu+\bfnu_x \, dx$,
$\bfnu_2=\bfnu+\bfnu_y \, dy$ the formulae (\ref{eq:6.13})
convert into the classical Lelieuvre formula, the relation
(\ref{eq:6.21}) is reduced to the first equation (\ref{eq:2.28}) and
\beq
\begin{array}{l}
det|{\bfnu},{\bfnu}_1,{\bfnu}_{2}| \rightarrow
det|{\bfnu},{\bfnu}_x,{\bfnu}_{y}| \, dx \, dy
\;, \\
det|{\bfnu},{\bfnu}_1,{\bfnu}_{11}|
\rightarrow
det|{\bfnu},{\bfnu}_x,{\bfnu}_{xx}| \, dx^3
\;, \\
det|{\bfnu},{\bfnu}_2,{\bfnu}_{22}|
\rightarrow
det|{\bfnu},{\bfnu}_y,{\bfnu}_{yy}| \, dy^3
\,\, .
\end{array}
\label{eq:6.22}
\ee
So we have
\begin{proposition}
The {\em l.h.s.} of (\ref{eq:6.22}) and (\ref{eq:6.17}) or better
\beq
\begin{array}{l}
\Omega_2 = 
\langle ( {\bff}_2-{\bff}) , ({\bfnu}_1 -
{\bfnu}) \rangle =
det|{\bfnu},{\bfnu}_1,{\bfnu}_{2}| \;, \\
\Omega_3 = 
\langle ( {\bff}_1-{\bff}_{-1}) , ({\bfnu} -
{\bfnu}_{-1}) \rangle =-
det|{\bfnu}_{-1},{\bfnu},{\bfnu}_{2}| \;, \\
\tilde{\Omega}_3 = 
\langle ( {\bff}_2-{\bff}_{-2}) , ({\bfnu} -
{\bfnu}_{-2}) \rangle =-
det|{\bfnu}_{-2},{\bfnu},{\bfnu}_{2}| \;\;. 
\end{array}
\label{eq:6.23}
\ee
are the discrete analogs of the Blaschke and Fubini cubic forms of affine 
surfaces.
\end{proposition}
Finally, we consider the dual affine gauge $\nu_4=1$.
In this case the inverse map (\ref{eq:5.1}) has a simple form
\beq
{\bfnu}_1 -{\bfnu} =\left[ {\bff},{\bff}_1 \right] \;\; ,\;\; 
{\bfnu}_2 -{\bfnu} =-\left[ {\bff},{\bff}_2\right]
\label{eq:6.24}
\ee
and instead of formulae (\ref{eq:6.14})-(\ref{eq:6.23})
one has those with the substitutions ${\bfnu} \leftrightarrow
{\bff}$.

A class of discrete surfaces for which both $f_4=-1$
and $\nu_4=1$, {\em i.e}
$\langle {\bff} , {\bfnu} \rangle=1$
is given by the discrete affine spheres. They have been studied 
recently in \bib{r10} where the formulae (\ref{eq:6.13})-(\ref{eq:6.16}),
(\ref{eq:6.24}) have been derived in this particular case. 

 \hfill


   \begin{centerline} 
   {\bf REFERENCES}
   \end{centerline}

   \begin{enumerate}

  \item \label{r1}
        Lelieuvre A., {\it Sur les lignes asymptotiques et leur
	representation spherique}, Bull. des Sciences Math. Astron. (2),
	{\bf 12}, (1988), 126-128.
  \item \label{r2}
	Blaschke W., {\it Differential Geometrie, II. Affine 
	differential geometrie}, Verlag, Berlin, 1923.
 \item \label{r3} 
	Schirokow P.A. and Schirokow A.P., {\it Affine
	differential geometrie}, Leipzig, 1962;
	{\it Affine differential geometrie}, Gosizdat., Moskow, 1959
	(in Russian).
  \item \label{r4}
	Nomizu K. and Sasaki T., {\it Affine differential geometry},
	Cambridge, Univ.Press, 1994.
  \item \label{r5}
	Li A.-M., Nomizu K. and  Wang C., {\it A generalization of
	Lelieuvre's formula}, Results in Math., {\bf 20}, (1991),
	682-690.
  \item \label{r6}
	Konopelchenko B.G. and Pinkall U., {\it Integrable deformations
	of affine surfaces via Nizhik-Veselov-Novikov equation},Physics Letters
        A (to appear),preprint SFB 288 No. 307, 1998.
  \item \label{r7}
	Fubini G. and Cech E., {\it Introduction a la geometrie projective
        differentielle des surfaces }, Gauthier-Villars, Paris, 1931.
  \item \label{r8}
	Bobenko A. and Pinkall U., {\it Discrete surfaces with constant
	negative Gaussian curvature and the Hirota equation}, J. 
	Differential Geometry, {\bf 43}, (1996), 527-611.
  \item \label{r9}
	Bobenko A.I. and Seiler R.(Eds.), Discrete integrable geometry
	and physics, Oxford Univ. Press, Oxford, 1998.
  \item \label{r10}
	Bobenko A. and Schief W., {\it Affine spheres:Discretization and
	Duality Relations}, preprint SFB 288 No. 297 (1997).

\end{enumerate}
\end{document}